\input amssym.def
\input amssym.tex

\magnification=\magstep1
\hsize=17,5truecm
\vsize=25.5truecm
\hoffset=-0.9truecm
\voffset=-0.8truecm
\topskip=1truecm
\footline={\tenrm\hfil\folio\hfil}
\raggedbottom
\abovedisplayskip=3mm 
\belowdisplayskip=3mm
\abovedisplayshortskip=0mm
\belowdisplayshortskip=2mm
\normalbaselineskip=12pt  
\normalbaselines

{\bf Fields of cohomological dimension one versus $C_1$-fields}

\bigskip

J.-L. Colliot-Th\'el\`ene

\bigskip

Summary. {\it  Ax gave examples of fields of cohomological
dimension 1 which are not $C_1$-fields. Kato and Kuzumaki
asked whether a weak form of the $C_1$-property holds
for all fields of cohomological dimension 1 (existence
of  solutions in extensions of coprime degree rather
than existence of a solution in the ground field).
Using work of Merkur'ev and Suslin, and of Rost,
D. Madore and I recently produced examples
 which show
that the answer is in the negative. In the present
note, I produce examples
  which require
less  work than the original ones. In the original paper,
some of the examples were given by forms of degree 3 in 4 variables.
Here, for an arbitrary prime $p \geq 5$, I use forms of degree $p$ in
$p+1$ variables.}

 \bigskip

{\bf Introduction}

\medskip

 Let $X$ be an algebraic variety over a field $k$. 
The index $I(X)$ of $X/k$ is the greatest common
denominator of the degrees over $k$ of the residue
fields at closed points of $X$ :  
$$I(X) = {\rm g. c. d.}_{x \in X_{0}} [k(x):k] .$$
This is also the greatest common divisor of the degrees
of finite field extensions $K/k$ such that the set $X(K)$
of $K$-rational points of $X$ is not empty.

If $X$ has a $k$-rational point, then $I(X)=1$,
but the converse does not generally hold.

\medskip

  A field $k$ is said to be $C_1$ if any 
homogeneous polynomial $F(x_0,\dots,x_n) \in k[x_0,\dots,x_n]$
 of degree $d$  in $n+1>d$ variables has a nontrivial zero, that is
there exist $\alpha_0, \dots, \alpha_n$ in $k$, not all
of them zero, such that $F(\alpha_0,\dots,\alpha_n)=0$.

The three basic examples of such fields are

1) Finite fields (Chevalley, Warning, 1935)

2) Function fields in one variable over an algebraically
closed field (Tsen, 1933)

3) Formal power series field in one variable over
an algebraically closed field (Lang, 1952)

\medskip

Let $k$ be a perfect field, $\overline k$ an algebraic
closure, $\frak{g}_k={\rm Gal}({\overline k}/k)$.
A perfect field $k$ is of (Galois) cohomological
dimension $\leq  1$ if any of the following
equivalent properties hold (see Serre [S]).

(i) For any (continuous) finite $\frak{g}_k$-module $M$ and
any integer $i \geq 2$, $H^{i}(\frak{g}_k,M)=0$.

(ii) For any finite field extension 
$K \subset {\overline k}$, the Brauer group
${\rm Br}(K)= H^2(\frak{g}_K,{\overline k}^*)$ vanishes.

(iii) Any homogeneous space $X/k$ under a connected linear
algebraic group  has a $k$-rational point.

If the group $\frak{g}_k$ is a pro-$p$-group, and $p \neq {\rm char.}(k)$,
 then these
conditions are equivalent to the mere condition

(iv) The $p$-torsion subgroup of ${\rm Br}(k)$
is trivial.

\medskip

If a (perfect) field $k$ is $C_1$, then it is of cohomological
dimension at most 1. In 1965, Ax [A] showed that the converse does
not hold. He produced an example of a field $k$
which is of cohomological dimension 1, and a form
$F$ of degree 5  in 10 variables over that field with no nontrivial zero.
However, the very construction of that form shows that it
possesses a zero in field extensions of degree 2, 3 and 5.
If we let $X$ be the hypersurface in 9-dimensional projective
space defined by this form, we have $I(X)=1$. Ax gave other
examples,
but they all have the property that the index
of the associated hypersurface is 1.

\medskip

In 1986, Kato and Kuzumaki [KK] asked : If a field $k$ is
of cohomological dimension at most~1, and 
$X \subset {\bf P}^n_k$ is a hypersurface defined
by a form of degree at most $n$, does it follow
that $I(X)=1$ ?

In other words, is the field $C_1$ as far as zero-cycles of
degree 1 are concerned (as opposed to rational points) ?

\medskip

In the  article [CM], David Madore and I showed that
the answer to this question is in the negative.
We produce  a field $k$ of cohomological dimension 1
and a cubic surface $X \subset {\bf P}^3_k$ 
such that $I(X)=3$. 
The geometric Picard
group of a smooth cubic surface is of rank 7.
This creates some technical difficulties.

In the present note, where I review the method
of [CM], I give new examples which are easier to
discuss. The result is the following :

\medskip

Theorem 1. {\it  For each prime $p \geq 5$, there exists
a field $F$ of characteristic zero, of cohomological
dimension 1 and a smooth hypersurface $X \subset {\bf P}_F^p$
of degree $p$ such that $I(X)=p$.}
\medskip

\bigskip

{\bf How to produce fields of cohomological dimension  $\leq 1$}

\medskip

Let $k$ be a field, $n\geq 2$ an integer. A Severi-Brauer variety
of index $n$ over $k$ is a twisted form of projective space
${\bf P}^{n-1}$, that is, it is a $k$-variety which after
a suitable extension $K/k$  becomes
isomorphic to ${\bf P}^{n-1}_K$.

There is a bijection, due to F. Ch\^atelet (1944),
 between the set of $k$-isomorphism
classes of such $k$-varieties and the set of isomorphism
classes of central simple $k$-algebras of index $n$.
In this bijection, projective space over $k$
corresponds to the matrix algebra. This is referred
to as the trivial class.
For $n=2$, this is the well-known correspondence
between conics and quaternion algebras.

\medskip

Theorem 2. {\it Let  $k$ be a field. Let $p$ be a
prime which does not divide the characteristic of $k$.
If one starts from $k$ and one
 iterates the following two operations

1) go from a field $K$ to the fixed field of
a pro-$p$-Sylow subgroup of the absolute Galois
group of $K$,

2) go from a field $K$ to the function field
of a nontrivial Severi-Brauer variety over $K$ of index $p$,

 then one ultimately obtains  a field $F$ containing
$k$ whose cohomological dimension is at most 1.}

\medskip

Sketch of proof : For the  field $F$ one obtains in
the limit, the operations in 1) ensure that the
absolute Galois group of $F$ is a pro-$p$-group.
To show that $F$ is of cohomological dimension 1,
it thus suffices to show that the $p$-torsion of
the Brauer group of $F$ is trivial. The Merkur'ev-Suslin
theorem (1983) implies that for any field $K$
of characteristic different from $p$ containing
the $p$-th roots of 1, the $p$-torsion of the
Brauer group is generated by the classes of
central simple algebras of index $p$. But
the operations in 2) ensure that over the field $F$
there is no such nontrivial central simple algebra.

\bigskip

{\bf  Galois action and the Picard group}

\medskip

Proposition 3. {\it Let $k$ be a field, ${\overline k}$
a separable closure of $k$, $\frak{g}={\rm Gal}({\overline k}/k)$.
Let $X/k$ be smooth, projective, geometrically integral 
variety. Write  ${\overline X}= X \times_k{\overline k}$.
Let $k(X)$ be the function field of $X$.
One then has a natural exact sequence
$$ 0 \to {\rm Pic}(X) \to {\rm Pic}({\overline X})^{\frak{g}}
\to {\rm Br}(k) \to {\rm Br}(k(X)).$$
The sequence is functorial contravariant with respect
to dominant $k$-morphisms.
}

\medskip

Here ${\rm Pic}(X)$ is the Picard group of $X$.
The second and fourth maps are the obvious
restriction maps. 

When $X$ is a Severi-Brauer variety corresponding
to a central simple 
algebra $A$ of index $m$
 then $ {\rm Pic}({\overline X})= {\rm Pic}({\bf P}^{m-1})=
{\bf Z} {\cal O}_{{\bf P}^{m-1}}(1) $
and the image of ${\cal O}_{{\bf P}^{m-1}}(1)$ is the class of $A$
 in the Brauer group of $k$. The kernel of the restriction
map ${\rm Br}(k) \to {\rm Br}(k(X))$ is the finite cyclic group
spanned by the class of $A$.

\medskip

 Exercise. Use the above sequence and its functoriality
to establish the following frequently rediscovered result. 

\medskip

Proposition 4. {\it  Let $C/k$ be a smooth conic. Let $f : C \to D$
be a $k$-morphism of smooth projective geometrically integral
curves. If the degree of $f$ is even, then $D(k)\neq \emptyset$.}

\medskip

(If $f$ is constant, the result is clear. If $f$ is not constant,
then by L\"uroth's theorem $D$ is of genus zero, hence is
a smooth conic.)

\medskip

Corollary 5. {\it Let $X \subset {\bf P}^n_k$ be a 
smooth hypersurface.
If $n \geq 4$, then the restriction map ${\rm Br}(k) \to {\rm Br} (k(X))$
is one-to-one. }

\medskip

Indeed, it is a theorem of Max Noether that
under the above assumptions the group ${\rm Pic}({\overline X})$
is free of rank one, spanned by the class of a hyperplane section.
Since such hyperplanes are defined over $k$, the result follows from
the above proposition.

\medskip

Remark. For surfaces in ${\bf P}^3$, the situation is more
complicated. This accounts for the more elaborate arguments
used in [CM] to produce examples with cubic surfaces.

\bigskip

{\bf Rost's degree formula}

\medskip

To any prime $p$ and any 
projective irreducible variety $X$ over a field $k$,
Rost associates a class $\eta_p(X) \in {\bf Z}/I(X)$.
This class is killed by $p$. If $X$ is a nontrivial 
 Severi-Brauer variety
dimension $p-1$, then $I(X)=p$ and $\eta_p(X)=1 \in 
{\bf Z}/p$.

The construction of this invariant belongs to the
world of coherent modules and intersection theory.
There is no Galois cohomology here.

\medskip

Theorem 6 (Rost, cf. Merkur'ev). {\it
Let $f : Z \to X$ be a $k$-morphism of
proper integral $k$-varieties of the
same dimension. Then $I(X)$ divides $I(Z)$
and $$\eta_p(Z) = {\rm deg}(f) \eta_p(X) \in {\bf Z}/I(X).$$
When $Z$ is moreover smooth, the same result
holds under the mere assumption that $f$ is
a rational map from $Z$ to $X$.}

\medskip

This formula implies in particular that
$\eta_p(X)=0$ for any
variety which can be written as a product
$X= Y \times Z$ where $Z$ is a $k$-variety
of dimension at least one and $I(Z)=1$.

\medskip

Exercise. Use Rost's degree formula to give an alternate
proof of Proposition 4.
This formula actually yields the following  more general result. 

\medskip
Proposition 7. {\it 
If $p$ is a prime and $f : Z \to X$ is a dominant rational map
from a Severi-Brauer variety $Z/k$ of dimension $p-1$ to
a projective integral $k$-variety $X$, and
$p$ divides the degree of $f$, then $I(X)=1$.
This is in particular so if the dimension of $X$
is less than that of $Z$.
}

\vfill\eject

{\bf The example}

\medskip

We can now prove Theorem 1. More precisely, we shall 
establish the following result.

Theorem 8. {\it Let $p \geq 5$ and l be distinct primes
such that $p$ divides $(l-1)$, that is ${\bf F}_l^*/{\bf F}_l^{*p}
\neq 1$. Let $\alpha \in {\bf Z}$ such that the class
of ${\alpha}$ in ${\bf F}_l$ is not a $p$-th power.
Let $X \subset {\bf P}^p_{\bf Q}$ be the smooth hypersurface
over ${\bf Q}$ defined by the equation
$$ x_1^p + l x_2^p + \dots + l^{p-1}x_p^p - \alpha x_0^p = 0.$$
There exists a field $F$ of characteristic zero, of cohomological
dimension 1, such that the index $I(X_F)$ of the $F$-variety
$X_F=X\times_{\bf Q}F$ is equal to $p$.}

\bigskip

Proof.
It is clear that $I(X)=I(X/{\bf Q})$ divides $p$.
Let $K/{\bf Q}_l$ be an extension 
 of the $l$-adic field ${\bf Q}_l$ of degree
prime to $p$.
One easily checks that $X(K)=\emptyset$.
Thus $p=I(X_{{\bf Q}_l})$, hence also $p=I(X)$.

Starting from $k={\bf Q}$ (or if one wishes
$k={\bf Q}_l$), one then applies the process
described
 in Theorem 2. To achieve the announced
result, one must check that under each of
the changes of fields described in that theorem,
the condition $p=I(X)$ is preserved.
This is obvious for the change 1), which consists
in going over to the fixed field of a pro-$p$-Sylow
subgroup. 

Let   $p=I(X/K)$, and let $Y/K$ is a nontrivial
Severi-Brauer variety of dimension $p-1$.
Assume $I(X_{K(Y)})=1$. Then there exist
a projective integral $K$-variety $Z$ of dimension
$p-1$, a dominant $K$-morphism $f :Z \to Y$ of degree
prime to $p$ and a morphism
 $h : Z \to X$.
Then $I(X)$ divides $I(Z)$ and $p=I(Y)$ divides $I(Z)$.
By Rost's degree formula we have
$$\eta_p(Z) = {\rm deg}(f) \eta_p(Y) \in {\bf Z}/I(Y)={\bf Z}/p,$$
and
$$\eta_p(Z) = {\rm deg}(h) \eta_p(X) \in {\bf Z}/I(X)={\bf Z}/p.$$
Since $\eta_p(Y)=1 \in {\bf Z}/p$ and the degree of $f$ is prime to
$p$, the first equality implies $$\eta_p(Z) \neq 0 \in {\bf Z}/p.$$
The second equality then implies that $p$ does not
divide the degree of $h$.
The restriction map of $p$-torsion groups
${}_p{\rm Br}(K) \to {}_p{\rm Br} (K(Z))$ factorizes 
as ${}_p{\rm Br}(K) \to {}_p{\rm Br}(K(X)) \to {}_p{\rm Br}(K(Z)).$
The first map is injective by Corollary 5. Since the degree
of $h$ is prime to $p$, a corestriction-restriction argument
shows that the second map is also an injection.
On the other hand the factorization
${}_p{\rm Br}(K) \to {}_p{\rm Br} (K(Y))  \to {}_p{\rm Br} (K(Z))$
shows that the class $A_Y \neq 0$ of the Severi-Brauer variety $Y$
in ${\rm Br}(K)$ vanishes in ${\rm Br}(K(Z))$, since this class
vanishes in ${\rm Br}(K(Y))$. This contradiction shows that
$I(X_{K(Y)})=p$, and this completes the proof that for
the field $F$ of cohomological dimension at most 1 which one
obtains in the limit, one has $p=I(X_F)$.

\bigskip

I refer to [CM] for more detailed literature references and  for comments
on  the general context in which such problems arise.

\vfill\eject

{\bf References}

\bigskip

[A]  J. Ax, A field of cohomological dimension 1 which is
not $C_1$, Bull. Amer.Math. Soc. {\bf 71} (1975), 717.

[CM] J.-L. Colliot-Th\'el\`ene et D. Madore, Surfaces de
del 
Pezzo sans point rationnel sur un corps de dimension
cohomologique 1, Journal de
l'Institut Math\'ematique de Jussieu (2004) {\bf 3} (1), 1-16.

[KK] K. Kato and T. Kuzumaki, The dimension of fields and
algebraic K-theory, J. Number Theory {\bf 24} (1986) 229-244.

[M] A. S. Merkurjev, Rost's degree
formula, http://www.math.ucla.edu/\~{}merkurev/

[S] J-P. Serre, Cohomologie galoisienne, cinqui\`eme \'edition,
r\'evis\'ee et compl\'et\'ee, Springer Lecture Notes in Mathematics, vol.
{\bf  5}  (Springer, 1994).

\vskip1cm

J.-L. Colliot-Th\'el\`ene

CNRS

UMR 8628

Math\'ematiques

B\^atiment 425

Universit\'e Paris-Sud

F-91405 Orsay

France

\medskip

e-mail : colliot@math.u-psud.fr

 \bye